\documentclass[10pt]{amsart}
\usepackage{amsmath}
\usepackage{amsthm}
\usepackage{amsfonts}
\usepackage{amssymb}
\usepackage[all]{xy}

\newcommand{\supp}{\text{supp}}
\newcommand{\col}{\colon}
\newcommand{\Sym}{\text{Sym}^m_{\Ps^5}\mathcal C}
\newcommand{\Symg}{\text{Sym}^{m_g}_{\Ps^5}\mathcal C}
\newcommand{\SymP}{\text{Sym}^{h_g}(\Ps^2)^\vee}
\newcommand{\B}{\overline{B_{m_g}}}
\newcommand{\Ph}{\mathbb P_{h_g}}
\newcommand{\Ps}{\mathbb{P}}
\newcommand{\co}{\mathbb C}
\newcommand{\ra}{\rightarrow}
\newcommand{\spa}{\text{span}}
\newcommand{\M}{\overline{M_{0,m_g}}}
\newcommand{\N}{\overline{N_{0,m_g}}}
\newcommand{\lra}{\longrightarrow}
\newcommand{\J}{\overline{J_g}}
\renewcommand{\phi}{\varphi}

    \newtheorem{Lem}{Lemma}[section]
    
    \newtheorem{Prop}[Lem]{Proposition}
    \newtheorem{Thm}[Lem]{Theorem}

   \theoremstyle{definition}
   \newtheorem{Def}[Lem]{Definition}
   \newtheorem{Exa}[Lem]{Example}

\begin{document}

\title{Compactifying moduli of hyperelliptic curves}

\author{Marco Pacini}

\address{Dipartimento di Matematica Guido Castelnuovo, Universit\`a Roma La Sapienza, piazzale Aldo Moro 2, 00185 Roma, Italia}

\address{Instituto de Matematica Pura e Aplicada, Estrada D. Castorina 110, 22460-320 Rio de Janeiro, Brazil}

\email{pacini@impa.br}

\begin{abstract}
We construct a new compactification of the moduli space $H_g$ of smooth hyperelliptic curves of genus $g.$ We compare our compactification with other well-known remarkable compactifications of $H_g$.
\end{abstract}

\maketitle

\section{Introduction}

Let $H_g$ be the moduli space of smooth hyperelliptic curves of genus $g\ge 3.$ Several compactifications of $H_g$ have been constructed. For example, there exists a moduli space $\overline{B_m}$ of GIT-semistable binary forms of degree $m$, where a binary form of degree $m$ is a homogeneous polynomial of degree $m$ in two variables over $\co$, up to non-trivial constants. 
In particular, $\overline{B_{2g+2}}$ contains $H_g$ as a dense open subset.  
Recall that a \emph{$m$-pointed stable curve of genus zero} $(Y,p_1,\dots,p_m)$ is a curve $Y$ of genus $0$ with an ordered set of distinct smooth points $p_i\in Y$ such that $|Y_i\cap\overline{Y-Y_i}|+|Y_i\cap\{p_1,\dots,p_m\}|\ge 3$ for every irreducible component $Y_i$ of $Y$. A \emph{$m$-marked stable curve of genus zero} is a $m$-pointed stable curve of genus zero $(Y,p_1,\dots,p_m)$ up to the action of the symmetric group $S_m$ on $p_1,\dots,p_m$. A natural compactification of $H_g$ is given by its closure within the moduli space of Deligne-Mumford stable curves. This compactification is isomorphic to the moduli space $\overline{N_{0,2g+2}}$ of $2g+2$-marked stable curves of genus zero. There exists also a fine moduli space $\overline{M_{0,m}}$ of $m$-pointed stable curves of genus zero and $\overline{N_{0,m}}=\overline{M_{0,m}}/S_m.$

As pointed out in \cite{AL}, $\overline{N_{0,2g+2}}$ and $\overline{B_{2g+2}}$ are different schemes. We construct a compactification $\J$ of $H_g$, given in term of configurations of plane lines and we compare it to $\overline{N_{0,2g+2}}$ and $\overline{B_{2g+2}}$. Indeed, consider $(C,p_1,\dots,p_{2g+2})$, $C$ a smooth plane conic and $p_i\in C$ distinct points. Pick the $h_g=\binom{2g+2}{2}$ lines spanned by $p_i,p_j$ for $1\le i< j\le 2g+2.$ By taking the closure within $Sym^{h_g}(\Ps^2)^\vee,$ we obtain configurations associated to 
$(C,p_1,\dots,p_{2g+2}),$ $C$ singular or $p_i=p_j$. The variety $\J$ is the GIT-quotient of the set of GIT-semistable configurations of lines, with respect to the action of $SL(3)$. A boundary point of $\J$ is a configuration containing at least a non-reduced line. For example, if $C$ is smooth and $p_1=p_2\ne p_j$ for $j\ge 3,$ the associated configuration contains $\spa(p_1,p_j)_{j\ge 3}$ as double lines. The boundary points of $\J$ have the following geometric meaning. If $C$ is smooth and $p_i\in C$ are distinct, consider the double cover $\phi:X\ra C$ branched at the $2g+2$ points $p_i.$ From \cite[pag. 288]{ACGH} and \cite[Proposition 6.1]{M}, we have that $\mathcal O_X(\phi^*(p_i+p_j))$ is a $(g-1)$-th root of $\omega_X$. If $g=3$, they are the 28 odd theta characteristics of the hyperelliptic curve $X$. Thus $\J$ is a compactification of $H_g$ given in terms of limits of configurations of higher spin curves of order $g-1$, in the sense of \cite{CapCasCorn}.

In \cite{H} and \cite {AL}, the authors construct a geometrical meaningful morphism $F_g:\overline{N_{0,2g+2}}\ra\overline{B_{2g+2}}$. In Theorem \ref{first}, Theorem \ref{second} and Theorem \ref{fact}, we construct rational maps $\overline{N_{0,2g+2}}\stackrel{\beta_g}{\dashrightarrow}\overline{J_g}\stackrel{\alpha_g}{\dashrightarrow}\overline{B_{2g+2}}$, giving a factorization of $F_g$. The construction of $\alpha_g$ follows from Lemma \ref{injection}, proving that it is possible to recover $(C,p_1,\dots,p_{2g+2}),$ $C$ smooth conic and $p_i\in C,$ from its configuration of lines. In particular, Lemma \ref{injection} extends the results of \cite{CS1} and \cite{L}, stating that a smooth plane quartic can be recovered from its bitangents, to double conics. In fact, the stable reduction of a general one-parameter deformation of a double conic $C$ is an hyperelliptic curve, the double cover of $C$ branched at $8$ points. The limits of the bitangents give rise to the configuration of lines associated to $C$ and the 8 points. We point out that a different generalization of bitangents for any genus is theta hyperplanes,  used in \cite{CS2} and \cite{GS}.

In short, in Section 2 we show properties of twisters of curves.  In Section 3, we construct $\J$ and the map $\alpha_g.$ In Section 4, we construct the map 
$\beta_g$, showing that $\alpha_g$ and $\beta_g$ provide a factorization of $F_g$.

\subsection{Notation}
We work over $\co$. A \emph{family of curves} is a proper and flat morphism $f\col\mathcal C\ra B$ whose fibers are curves. If $0$ is a point of a scheme $B,$ set $B^*:=B-0$. A \emph{smoothing} of a curve $C$ is a family $f\col\mathcal C\ra B,$ where $B$ is a smooth, connected, affine curve of finite type, with a distinguished point $0\in B,$ such that $f^{-1}(0)$ is isomorphic to $C$ and $f^{-1}(b)$ is smooth for $b\in B^*.$ A \emph{general smoothing} is a smoothing with smooth total space. Let $Y$ be a scheme and $X\ra Y$ be a $Y$-scheme. Denote by $Sym^m_Y X:=X_Y^m/S_m$ the quotient of $X_Y^m=X\times_Y X\times_Y\dots\times_Y X$ (the $m$-fiber product) by the symmetric group $S_m$. If a group $G$ acts on a variety $X$, denote by $X^{ss}$ the set of GIT-semistable points. If $p,q\in\mathbb{P}^2$, $p\ne q$, we set $\overline{pq}=\spa(p, q)$. For a positive integer $g$, set $m_g=2g+2$ and $h_g=\binom{m_g}{2}$.

\section{On some properties of conic twisters}

A \emph{twister} of a curve $Y$ is a $T\in\text{Pic}(Y)$ such that there exists a smoothing $\mathcal Y$ of $Y$ such that $T\simeq\mathcal O_{\mathcal Y}(D)\otimes\mathcal O_Y,$ where $D$ is a Cartier divisor of $\mathcal Y$ supported on irreducible components of $Y.$ If $Y$ is of compact type, it is well-known that a twister depends only on its multidegree. A \emph{conic twister} is a twister whose degrees on the irreducible components are positive and sum up to 2.

\begin{Prop}\label{conictwister}
Let $Y$ be a genus zero curve and $T$ a conic twister of $Y$. Then:
\begin{itemize}\label{twister}
\item[(i)]
if $d_1,\dots,d_N$ are positive integers summing up to 2, then there exists a conic  twister $T$ of $Y$ such that $\underline{\deg}\;(\omega_Y^\vee\otimes T)=(d_1,\dots,d_N)$;
\item[(ii)]
the linear system $|\omega_Y^\vee\otimes T|$ is base point free, two-dimensional and induces a morphism $Y\ra\Ps^2$ realizing $Y$ as plane conic.
\end{itemize}
\end{Prop}

\begin{proof}

(i) Given two components $Y_1, Y_2\subset Y$ such that $Y_1\cap Y_2=\{p_1\sim p_2\}$, the class $[p_1]-[p_2]$ is a twister $T$ such that $T|_{Y_1}=\mathcal O_{Y_1}(1)$, $T|_{Y_2}=\mathcal O_{Y_2}(-1)$ and $T$ is trivial on the other components of $Y$; the claim follows from the connectivity of $Y$.  

(ii) If $\deg_{Y_1}T=2$ for some $Y_1\subset Y$, then $|\mathcal O_{Y_1}(T|_{Y_1})|=|\mathcal O_{\Ps^1}(2)|\simeq\Ps^2$, and the map 
$Y_1\ra\Ps^2$ is degree 2. This map extends to $Y$ because the dual graph of the components of $Y$ is a tree. If $\deg_{Y_1} T=\deg_{Y_2} T=1$ for some $Y_1,Y_2\subset Y$, let  $\Gamma$ be the unique path in the dual graph of the components of $Y$ connecting $Y_1,Y_2$, and let $p_1\in Y_1$, $p_2\in Y_2$ be the unique points in $Y_1,Y_2$ which sit on the dual of $\Gamma$ in $Y$. Note that $T$ induces isomorphisms $\pi_i\col Y_i\overset{\simeq}{\ra}|\mathcal O_{Y_i}(T|_{Y_i})|\simeq\Ps^1$ for $i=1,2$. Since the dual graph of the components of $Y$ is a tree, there is a unique map $\pi\col Y\ra\Ps^1\cup_{\pi_1(p_1)\sim\pi_2(p_2)}\Ps^1$ extending $\pi_1,\pi_2$. Since up to projective morphisms there is a unique embedding $\Ps^1\cup_{\pi_1(p_1)\sim\pi_2(p_2)}\Ps^1\ra\Ps^2$, we are done.
\end{proof}

\section{The first map}

Let $\mathcal C\ra\Ps^5\simeq|\mathcal O_{\Ps^2}(2)|$ be the universal plane conic. For any integer $m\ge 2$, consider the variety $\Sym$ and the morphism 
$\rho:\Sym\ra\Ps^5$. If $k\in\Sym$, let $\supp(k)$ be the conic parametrized by $\rho(k)$. The points of $k$ are called \emph{markings} and $(k,\supp(k))$, a \emph{conic with markings}. A marking has a \emph{weight}, i.e. the number of times it appears in $k$. We call the markings $p_{min}$ and $p_{max}$ with minimal and maximal weight, the \emph{minimal} and \emph{maximal markings} of $k$.  
Recall that $m_g=2g+2$ and $h_g=\binom{m_g}{2}$, where $g\ge 3$. 
Set $\Ph=\SymP$. Consider the rational map: 
\begin{equation}\label{psi1}
\psi:\Symg\dashrightarrow\Ph
\end{equation}

where, if $k$ has markings $\{p_i\}_{1\le i\le s}$ of weight $m_i=1$ and $\{p_i\}_{s<i\le r}$ of weights $m_i>1$, then:
\begin{equation}\label{psi2}
\psi(k)=(\cdots\underset{m_i\cdot m_j \text{ times}}{\underbrace{\overline{p_i p_j},\dots\overline{p_i p_j}}},\cdots\underset{\binom{m_h}{2} \text{ times}}{\underbrace{\mathbb{T}_{p_h}\supp\;k,\cdots\mathbb{T}_{p_h}\supp\;k}},\cdots)_{\underset{s<h\le r}{1\le i<j\le r}}.
\end{equation}

Let $\Gamma_{\psi}$ be the closure in $\Symg\times\Ph$ of the graph of $\psi$ and $p\col\Gamma_{\psi}\ra\Ph$ be the second projection. Consider the 
GIT quotient: $$q:p(\Gamma_{\psi})\cap\Ph^{ss}\lra\J=(p(\Gamma_{\psi})\cap\Ph^{ss})/SL(3).$$ 

We say that $k$ is \emph{degenerate} if $\text{supp(k)}$ is not integral. Consider the open subset $V$ of $p(\Gamma_\psi)\cap\Ph^{ss}$ defined as: $$V=\{r\in p(\Gamma_\psi)\cap \Ph^{ss}:r\ne p(k,r)\;\forall\;k \text{ degenerate}\}\subset p(\Gamma_\psi)\cap\Ph^{ss}$$

Recall the Hilbert-Mumford criterion \cite[Proposition 4.3]{MFK} for configurations of plane lines. Let $r$ be in $\Ph.$ For a point $p\in\Ps^2,$ let $\mu_p(r)$ be the number of lines of $r$, with multiplicities, containing $p.$ Let $\mu_l(r)$ be the multiplicity of a line $l$ of $r$. Then $r$ is GIT-semistable iff 
$\text{max}_{p\in\Ps^2}\,\mu_p(r)\le 2h_g/3$ and $\text{max}_{l\in(\Ps^2)^\vee}\mu_l(r)\le h_g/3$.

\begin{Lem}\label{unstable}
Let $(k, r)\in\Gamma_{\psi}$. Then:
\begin{itemize}
\item[(i)]
if $k$ has a marking $q$ of weight at least $g+1$, then $r$ is not GIT-semistable;
\item[(ii)]
if $\text{supp}(k)$ is reducible and the set of markings on smooth points of a component is one marking of weight 1, then $r$ is not GIT-semistable;
\item[(iii)]
if $\text{supp}(k)$ is integral and the markings have weight 1, then $r\in V$. 
\end{itemize}
\end{Lem}

\begin{proof}
(i) We have $\text{max}_{p\in\Ps^2}\,\mu_p(r)\ge\mu_q(r)\ge\binom{g+1}{2}+(g+1)^2>2h_g/3$.

(ii) From (i), we can assume that the node $n$ of $\text{supp}(k)$ is a marking of weight at most $g$. The number of lines of $r$ not containing $n$ is at most $2g+1$ and hence $\text{max}_{p\in\Ps^2}\,\mu_p(r)\ge\mu_n(r)\ge h_g-2g-1> 2h_g/3$. 

(iii) We have that $r$ is GIT-semistable, because $\text{max}_{p\in\Ps^2}\mu_p(r)=2g+1<2h_g/3$ and $\text{max}_{l\in(\Ps^2)^\vee}\mu_l(r)=1<h_g/3$. The property $\text{max}_{l\in(\Ps^2)^\vee}\mu_l(r)=1$ characterizes the configurations of integral conics with markings of weight 1, hence $r\in V$.
\end{proof}

\begin{Lem}\label{injection}
Consider the rational map $\psi:\Symg\dashrightarrow\Ph$. Then the restricted morphism $\psi:\psi^{-1}(V)\ra V$ is injective for every $g\ge 3$. 
\end{Lem}

\begin{proof} 
Pick $k\in\text{Sym}^m_{\Ps^5}\mathcal C$, where $m\ge 2$. Let $C=\supp\;k$ be integral. Set $r=\psi(k)$, as in (\ref{psi2}). Let $\{m_1,\dots,m_r\}$ be the set of the weights of $k$, where $m_i\le m_{i+1}$.

\smallskip

\noindent
{\bf Step 1.} 
Assume that $\{m_1,\dots,m_r\}\ne\{1,1\}$. The goal of the first step is to  recover the maximal markings of $k$ and their weights. We claim that the maximal markings of $k$ are the points $p\in\Ps^2$ with maximum multiplicity $\mu_p(r)$. It is easy if $m_i=1$, for $1\le i\le r$, thus assume that $m_{\text{max}}=\text{max}\;\{m_i\}_{1\le i\le r}\ge 2$. If $p\in C$, then $\mu_p(r)\le\mu_{p_{\text{max}}}(r)$, with the equality iff $p$ is a maximal marking of $k$ and we are done. If $p\notin C$, take two markings $p_i,p_j\in C$ of $k$ of weights $m_i$ and $m_j$ such that $p\in\overline{p_i p_j}$ if $p_i\ne p_j$ and $p\in\mathbb{T}_{p_i}C$ if $p_i=p_j$. Thus:

\begin{itemize}
\item[i)] 
if $p_i\ne p_j$ and $p_i,p_j\ne p_{\text{max}},$ then $\mu_{\overline{p_i p_j}}(r)=m_i m_j<m_{\text{max}}(m_i+m_j)=\mu_{\overline{p_i p_{\text{max}}}}(r)+\mu_{\overline{p_j p_{\text{max}}}}(r)$;
\item[ii)]
if $p_j=p_{\text{max}}$ and $p_i\ne p_{\text{max}},$ then $\mu_{\overline{p_i p_{\text{max}}}}(r)<\mu_{\overline{p_i p_{\text{max}}}}(r)+\mu_{\mathbb{T}_{p_{\text{max}}}C}(r)$;
\item[iii)]
if $p_i=p_j\ne p_{\text{max}},$ then $\mu_{\mathbb{T}_{p_i}C}(r)=\binom{m_i}{2}<m_i m_{\text{max}}=\mu_{\overline{p_i p_{\text{max}}}}(r)$;
\item[iv)]
two lines $\overline{p_i p_{\text{max}}}$ in two different cases among $i)$ $ii)$ $iii)$ cannot be the same;
\item[v)]
if $p_i=p_j=p_{\text{max}},$ then $\mathbb{T}_{p_{\text{max}}}C$ contains both $p_{\text{max}}$ and $p$. The case ii) does not hold.  
If $\mathbb{T}_{p_{\text{max}}}C$ is the only line of $r$ containing $p$, then $\mu_p(r)=\mu_{\mathbb{T}_{p_{\text{max}}}C}(r)<\mu_{\mathbb{T}_{p_{\text{max}}}C}(r)+\mu_{\overline{p_h p_{\text{max}}}}(r)\le\mu_{p_{\text{max}}}(r)$ for some $C\ni p_h\ne p_{\text{max}}$. If $p$ is in at least 2 lines of $r,$ then at least one case i) or iii) holds.
\end{itemize}

This shows that $\mu_p(r)<\mu_{p_{max}}(r)$. Thus, we recover the maximal markings of $k$ as the points of $\Ps^2$ with maximum multiplicity. In particular, we find also the number $N$ of maximal markings of $k$. To recover their weight $m_{max}$, consider $m=\sum_{1\le i\le r}m_i$ and the subconfiguration $r'$ of $r$ of the lines containing no maximal markings of $k$. If $r'\ne\emptyset$, then $r'$ is the configuration of lines associated to $(k',C)$, where $k'$ are the non-maximal weights of $k$. We know the sum of the multiplicities of the lines of $r'$, thus we know also the sum $m'$ of the weights of $k'$ and $m_{\text{max}}=(m-m')/N$. If $r'=\emptyset$, then either $m_i=m_j$ for $1\le i,j\le r$ and $m_{\text{max}}=m/N$, or $m_1=1<m_i=m_j$ for $1<i,j\le r$ and $m_{\text{max}}=(m-1)/N.$

\smallskip
\noindent
{\bf Step 2.} Pick $k\in\psi^{-1}(V)$ for $m=2g+2$. We recover the markings of $k$, with the exception of the marking of multiplicity $m_1$, if $m_1=1<m_2$, and of multiplicity $m_1$ and $m_2$, if $m_1=m_2=1<m_3$. In fact, using Step 1 we find the maximal markings of $k$ and their weights. Now, consider $(k',C)$, where $k'$ are the non-maximal markings of $k$. We find the maximal markings of $k'$ and their weights, using Step 1. By iterating, we find the configuration $r_0$ associated to the markings of $k$ with minimal weights. If either $m_1=m_2\ne 1$ or $m_1=m_2=m_3=1$, we find the markings of $k$ and their weights. Otherwise, let $\{p_1,\dots, p_s\}$ be the set of the recovered markings, where $s\ge 2$, by Lemma \ref{unstable} (i). Consider the subconfiguration $r''$ of $r$ obtained by getting rid of the lines containing $p_3,\dots, p_s$ and $\overline{p_1 p_2}$. We have three cases. In the first case, $r_0=\emptyset$ and $m_1=1<m_2$.  In the second case, $r_0$ is a line $l$ of multiplicity $\mu_l(r_0)>1$ and $3\le m_1<m_2$. The marking with weight $m_1$ is the point  contained in 3 lines of $r''$. We recover also its weight $m_1$, because $\mu_l(r_0)=\binom{m_1}{2}$. In the third case, $r_0$ is a line $l$ of multiplicity $\mu_l(r_0)=1$ and either $m_1=2<m_2$ or $m_1=m_2=1<m_3$. We have $m_1=2<m_2$ iff $r''$ has 5 lines. The marking with weight 2 is the unique point contained in 3 lines of $r''$. 

\smallskip
\noindent
{\bf Step 3.} Pick $k\in\psi^{-1}(V)$ for $m=2g+2$. If, using Step 2, we find at least 5 markings, we recover $(k, C)$, because we find also the markings with multiplicity 1 as the points on $C$ with multiplicity $2g+1$. Assume that, using Step 2, we find at most 4 markings and their weights sum up to $2g+2$, i.e. they are all the markings of $k$. They are at least 3 markings of weights at least 2, by Lemma \ref{unstable} (i) and Step 2. Their tangents are the lines which do not contain a pair of markings. We find at least 3 markings and 3 tangents to the markings, then we recover $(k,C)$. If, using Step 2, we find at most 4 markings of $k$ and their weights do not sum up to $2g+2$, then either $m_1=1<m_2$ or $m_1=m_2=1<m_3$. There are 5 cases. Notice that $k$ has at least 4 markings, by Lemma \ref{unstable} (i).

\smallskip
\noindent
a) 2 recovered markings and the weights are $\{1,1,m_3,m_4\}$, $2\le m_3\le m_4$.

Let $p_3,p_4$ be the markings with multiplicities $m_3,m_4$. If $m_3\ne 3,$ then $\mathbb T_{p_3} C$ is the line through $p_3$ not containing $p_4$ and whose multiplicity is not $m_3$. Similarly, we determine also $\mathbb T_{p_4} C$. Consider the 4 lines of the configuration $r$ different from $\mathbb T_{p_3} C,\mathbb T_{p_4} C, \overline{p_3p_4}$ and containing either $p_3$ or $p_4$. The pairwise intersections of these lines are 6 points: two points of multiplicity $2g+1$, two of multiplicity $2g$ and $p_3,p_4$. The points of multiplicity $2g+1$ are the markings with weight 1. Thus we recover 4 markings and 2 tangents to the markings, hence also $(k,\supp\;k)$.
If $m_3=3,$ then also $m_4=3$ from Lemma \ref{unstable} (i) and $g=3.$ It is easy to see that there is only one conic with markings having $r$ as associated configuration.

\smallskip
\noindent
b) 3 recovered markings and the weights are $\{1,m_2,m_3,m_4\}$, $2\le m_2\le m_3\le m_4$.

Let $p_i$ have weight $m_i$. Assume that two weights are not equal to 3, e.g.  $m_2,m_3\ne 3$. Consider the subconfiguration $r'$ of $r$ of the lines not containing $p_4$ and different from $\overline{p_2p_3}$.
The marking with weight 1 is the point $p\ne p_2,p_3$ with multiplicity $\mu_p(r')=m_2+m_3$.
If two weights are equal to 3, then $m_2=m_3=3$, $m_4=2g-5$. If $m_4\ne 3$,  then $\mathbb{T}_{p_4} C$ is the unique line containing $p_4$ with multiplicity 
$\binom{m_4}{2}$ and different from $\overline{p_2p_4},\overline{p_3p_4}$. Consider the subconfiguration $r'$ of $r$ obtained by getting rid of $\mathbb{T}_{p_4} C$.
Then $r'$ is the configuration of $(k', C)$, where $k'$ are the markings of $k$ with set of weights $\{1,1,3,3\}$, where $p_4$ has weight 1 in $k'$. 
Using (a), we recover $(k', C)$, hence also the original $(k,C)$. If $m_2=m_3=m_4=3$, it is easy to see that there is only one conic with markings having $r$ as associated configuration.

\smallskip
\noindent
c) 3 recovered markings and the weights are $\{1,1,m_3,m_4,m_5\}$, $2\le m_3\le m_4\le m_5$. 

There is a marking of weight different form 3, otherwise $2g+2=11$, for example $m_3\ne 3$.  Let $r'$ be the subconfiguration of $r$ obtained by getting rid of lines containing the points with weights $m_4,m_5$. Thus the markings with weight 1 are the points $p,q$ of multiplicity $\mu_p(r')=\mu_q(r')=m_3+1$. 
We recover $(k, C)$.

\smallskip
\noindent
d) 4 recovered markings and the weights are $\{1,m_2,m_3,m_4,m_5\}$, $2\le m_2\dots\le m_5$. 

Consider the subconfiguration $r'$ of $r$ obtained by getting rid of the lines connecting any two points with multiplicity $>1$. Then the point with multiplicity 1 is the only point contained in 4 lines of $r'$. We recover $(k, C)$.

\smallskip
\noindent
e) 4 recovered markings and the weights are $\{1,1,m_3,\dots,m_6\}$, 
$2\le m_3\dots\le m_6$.  

If a marking has a weight different form 3, we argue as in c). Otherwise, consider the subconfiguration $r'$ of $r$ obtained by getting rid of the lines containing the markings with weight $m_5,m_6$. Thus $r'$ is the configuration of $(k', C)$, where $k'$ has weights $\{1,1,3,3\}$ and we argue as in a).
\end{proof}

Denote by $J_g$ the subset of $\J$ corresponding to the classes of configurations of lines of integral conics with markings of weight 1.

\begin{Thm}\label{first} 
There exists a map $\alpha_g\col\J\dashrightarrow\B$ defined at least over $J_g$. 
\end{Thm}

\begin{proof}
Up to restrict $V$, we have that $\psi:\psi^{-1}(V)\ra V$ is an isomorphism by Theorem \ref{injection}. Pick $k\in \psi^{-1}(V)$. Since $\supp(k)$ is irreducible, Lemma \ref{unstable} (i) implies that $(k,\supp(k))$ is a GIT-stable binary form. Thus $\mathcal U\ra\psi^{-1}(V)\simeq V$ is a family of stable binary forms, where $\mathcal U=\{(k,p):p\in\supp\;k\}\subset \psi^{-1}(V)\times\Ps^2$, hence we get a $SL(3)$-invariant morphism $V\ra\B$, inducing the rational map 
$\alpha_g\col\J\dashrightarrow\B$. 

To show that $\alpha_g$ is defined over $J_g$, it is enough to show that the differential of $\psi$ is injective over an irreducible conics with markings of weight 1. We show that, if $U$ is the open subset of $\text{Sym}^{m_g}(\Ps^2)$ of $m_g$ distinct points such that any three of them are not contained in a line, then the differential of $\psi\col U\ra\Ph$ is injective, where $\psi(p_1,\dots p_{m_g})=(\dots\overline{p_i p_j}\dots)_{1\le i<j\le m_g}$. If $k\in U$, set $X=\psi(k)$, the union of $h_g$ distinct lines. If $\mathcal N_{X/\Ps^2}$ is the normal sheaf of $X$ 
in $\Ps^2$, then: $$\mathbb{T}_k U=\mathbb{T}_{p_1}\Ps^2\oplus\dots\oplus\mathbb{T}_{p_{m_g}}\Ps^2\stackrel{d\psi}{\lra}H^0(X,\mathcal N_{X/\Ps^2}).$$ For $v_i\in\mathbb{T}_{p_i}\Ps^2$, let $d\psi(v_1,\dots v_{m_g})=0$, i.e. it is the trivial embedded deformation of $X$, fixing all the components of $X$. This means that $v_i$ is contained in the lines of $X$ containing $p_i$. It is impossible, being $p_i$ contained in at least two lines of $X$.
\end{proof}

\section{The second map and the factorization}

A family of $m$-pointed stable curves of genus zero is a family $f:\mathcal Y\ra B$ of curves of genus zero with sections $\sigma_1,\dots,\sigma_m$ of $f$ such that $(Y_b,\sigma_1(b),\dots,\sigma_m(b))$ is a $m$-pointed stable curve of genus zero for $b\in B$. If $T$ is a conic twister of $Y$, let $\phi_T:Y\ra\Ps^2$ be the morphism induced by $|\omega_Y^\vee\otimes T|$ as in Lemma \ref{twister} (ii). 

\begin{Def}\label{princpart}
Let $(Y,p_1,\dots,p_m)$ be a $m$-pointed stable curve of genus zero. A connected subcurve $P\subset Y$ is \emph{a principal part} if there exists a conic twister $T$ of $Y$ such that the point $k=(\phi_T(p_1),\dots,\phi_T(p_m)\in\Sym$ 
satisfies $\supp\;k=\phi_T(P)=\phi_T(Y)$ and $\psi(k)\in\Ph^{ss}$, where 
$\psi:\Symg\dashrightarrow\Ph$ is the map (\ref{psi1}).
\end{Def}

A principal part $P$ has at most two components. If $P$ is a principal part of $(Y,p_1,\dots,p_{m_g}),$ the associated conic twister $T$ is uniquely determined by the condition $\phi_T(P)=\phi_T(Y).$ Since $P$ is connected, $\phi_T(p_i)$ is not in the singular locus of $\phi_T(P)$ and hence  
the map $\psi$ is defined over $k=(\phi_T(p_1),\dots,\phi_T(p_m))$.

\begin{Exa}\label{basic-exa}
Let $(Y,p_1,\dots,p_{10})$ be a pointed stable curve where $Y=Y_1\cup Y_2$. Assume that $p_1,\dots,p_6\in Y_1$ and $p_7,\dots,p_{10}\in Y_2$.
Both $Y$ and $Y_1$ are principal parts. In fact, if we consider $k_i=(\phi_{T_i}(p_1),\dots,\phi_{T_i}(p_{10}))$, where $T_1=\mathcal O_Y$ and $T_2$ is the twister given by $Y_2$, then $\psi(k_1)$ and $\psi(k_2)$ are GIT-semistable.  
\end{Exa}

We refer to \cite{K} for a proof of the following Lemma.

\begin{Lem}\label{Knudsen}
Let $[f:\mathcal Y\ra B,\sigma_1,\dots,\sigma_m]$ be a family of $m$-pointed stable curves of genus zero. Then there exists a unique family $[f':\mathcal Y'\ra B,\sigma'_1,\dots,\sigma'_{m-1}]$ and a $B$-morphism $h\col\mathcal Y\ra\mathcal Y'$ such that $h\circ\sigma_i=\sigma'_i$ for $i=1,\dots,m-1$. If $E_b\subset f^{-1}(b)$ is the component with $\sigma_m(b)\in E_b$ and $h_b=h|_{f^{-1}(b)}:f^{-1}(b)\ra (f')^{-1}(b)$, then:  
\begin{itemize}
\item[(i)]
$h_b$ contracts $E_b$ iff $\left|E_b\cap \overline{(f^{-1}(b)-E_b})|+|E_b\cap\{\sigma_1(b),\dots,\sigma_{m-1}(b)\}\right|\le 2$;
\item[(ii)]
if $h_b$ does not contract $E_b$, then $h_b$ is an isomorphism;
\item[(iii)]
if $h_b$ contracts $E_b$, then $h_b|_{\overline{f^{-1}(b)-E_b}}$ is an isomorphism.
\end{itemize}
\end{Lem}

\begin{Lem}\label{factorization}
The subset of $\M$ of the curves with a principal part is an open subset containing the locus of the curves with at most two components.
\end{Lem}

\begin{proof}
Let $P$ be a principal part of $(Y, p_1,\dots,p_{m_g})$. 
Let $[f:\mathcal Y\ra B,\sigma_1,\dots,\sigma_{m_g}]$ be a family of 
$m_g$-pointed stable curves of genus zero with $Y=f^{-1}(0)$ and $p_i=\sigma_i(0)$ for $0\in B$. Applying Lemma \ref{Knudsen}, we get a family $f':\mathcal Y'\ra B$ with $P=(f')^{-1}(0)$ and a morphism $h\col\mathcal Y\ra\mathcal Y'.$ Now, $P$ has at most 2 components, then up to  shrinking $B$ to an open subset containing $0$, all the fibers of $\mathcal Y'$ have at most two components. 
The image of the map $\phi:\mathcal Y'\ra\Ps(H^0(f'_*\omega_{f'}^\vee)^\vee)$ induced by $|\omega^\vee_{f'}|$ is a family of conics over $B$. The conic over 
$b\in B$ has markings $\gamma(b)=(\gamma_1(b),\dots,\gamma_{m_g}(b))$, where $\gamma_i=\phi\circ h\circ\sigma_i$. By construction a marking is not a node of a fiber, hence $\psi:\Symg\dashrightarrow\Ph$ is defined over $\gamma(b)$. Now, $\psi\circ\gamma(0)\in \Ph^{ss}$ because $P$ is a principal part. Thus, up to shrinking again $B$, we have $\psi\circ\gamma(b)\in\Ph^{ss}$ for $b\in B$ and the fiber of $\mathcal Y'\ra B$ over $b$ is a principal part of $f^{-1}(b)$ and we are done.

We show that $Y$ has a principal part if it has at most two components. If $Y$ is irreducible and $T=\mathcal O_Y$, 
then $\psi(\phi_T(p_1),\dots, \phi_T(p_{m_g}))$ is GIT-stable by Lemma \ref{unstable} (iii), thus $Y$ is a  principal part. If $Y=Y_1\cup Y_2$, let $[f:\mathcal Y\ra B, \sigma_1,\dots,\sigma_{m_g}]$ be a general smoothing of $(Y, p_1,\dots, p_{m_g})$ and $\mathcal C^*\ra B^*$ be the family of conics given by 
$|(\omega_f|_{f^{-1}(B^*)})^\vee|$. We have a map $\phi\col\mathcal Y\dashrightarrow\mathcal C^*$, which is an isomorphism away from $Y$. 
Now, $\psi(\phi\circ\sigma_1(b),\dots,\phi\circ\sigma_{m_g}(b))\in\Ph^{ss}$ for $b\in B^*$ by Lemma \ref{unstable} (iii). By the GIT-semistable replacement property, up to a finite base change totally ramified over $0\in B$, we find a completion $f'\col\mathcal C\ra B$ and $k\in\text{Sym}^{m_g}(C)$, where $C=(f')^{-1}(0)$ such that $\psi(k)\in\Ph^{ss}$.  If $\phi$ induces a morphism $\phi\col\mathcal Y\ra\mathcal C$, then $\phi|_Y=\phi_T$, where either $T=\mathcal O_Y$ or $T$ is given by $Y_1$ or $Y_2$ and $Y$ has a principal part.  Otherwise, call $\widetilde{\phi}:\mathcal{\widetilde{Y}}\ra\mathcal C$ the regularization of $\phi$. Let $\phi$ be not defined at $p\in Y$ and $E$ be an exceptional components over $p$. Let $p$ be a smooth point of $Y$. Then $\widetilde{\phi}(E)\subset C$ is a component containing at most one marking with weight 1 and $\psi(k)\notin\Ph^{ss}$ by Lemma \ref{unstable} (ii), a contradiction. Assume that $p_1,\dots,p_t\in Y_1$ for $t\le m_g/2$. Let $p$ be the node of $Y$. If $\widetilde{\phi}$ contracts $Y_2$, then $k$ has a marking with weight at least $m_g/2$ and $\psi(k)\notin\Ph^{ss}$ by Lemma \ref{unstable} (ii), a contradiction. If $Y_2$ is not contracted, then $C=\widetilde{\phi}(Y_2\cup E)$. If $q=\widetilde{\phi}(p_1)=\dots=\widetilde{\phi}(p_t)$, then $k=(q, \widetilde{\phi}(p_{t+1}),\dots,\widetilde{\phi}(p_{m_g}))$ and $\psi(k)\in\Ph^{ss}$. 
Consider $T=\mathcal O_T$ and $k'=(\phi_T(p_1),\dots,\phi_T(p_{m_g}))$. We have $\psi(k')\in\Ph^{ss}$ because $\text{max}_{p\in\Ps^2}\,\mu_p(\psi(k'))\le\text{max}_{p\in\Ps^2}\,\mu_p(\psi(k))$ and $\text{max}_{l\in(\Ps^2)^\vee}\mu_l(\psi(k'))\le\text{max}_{l\in(\Ps^2)^\vee}\mu_l(\psi(k))$, hence $Y$ is a principal part.
\end{proof}

Let $P_g\subset\M$ be the open subset of the curves with a principal part and 
$\N=\M/S_{m_g}$ the moduli space of $m_g$-marked stable curves, 
where $S_{m_g}$ is the symmetric group.

\begin{Thm}\label{second}
There exists a map $\beta_g\col \N\dashrightarrow \J$ defined at least over $P_g/S_{m_g}$. 
\end{Thm}

\begin{proof}
First of all, we show that if $(Y,p_1,\dots,p_{m_g})\in\M$ has two principal parts $P_1$ and $P_2$, whose associated configurations of lines are $r_1$ and $r_2$, then: 
\begin{equation}\label{orbit}
\overline{O_{SL(3)}(r_1)}\cap\overline{O_{SL(3)}(r_2)}\cap \Ph^{ss}\ne\emptyset,
\end{equation}

where $O_{SL(3)}(\cdot)$ denotes the orbit under the action of $SL(3)$. In fact, consider a general smoothing $[f:\mathcal Y\ra B,\sigma_1,\dots,\sigma_{m_g}]$ of $(Y,p_1,\dots,p_{m_g})$. For $j=1,2$, let $T_j$ be the twister of $P_j$. The morphisms $\phi_{T_j}$ induced by $|\omega_f^\vee\otimes T_j|$ give rise to two families of conics, which are isomorphic away from the special fiber. Furthermore, $\{\phi_{T_1}\circ\sigma_i\}$ and $\{\phi_{T_2}\circ\sigma_i\}$ induce markings on the two families. By construction, the associated families of configurations of lines are $SL(3)$-conjugate over $B^*$ and $r_1$ and $r_2$ are their  special fibers. Thus, $r_1$ and $r_2$ are GIT-semistable limits of conjugate families of GIT-stable configurations, hence (\ref{orbit}) follows.

Now, let $\mathcal P_g$ be the functor of familes of $m_g$-pointed stable curves with a principal part. We construct a functor tranformation $\mathcal P_g\ra\mathcal{M}or(-,\J)$. For a scheme $B$, pick $[f:\mathcal Y\ra B,\sigma_1,\dots,\sigma_{m_g}]\in\mathcal P_g(B)$. As in the proof of Lemma \ref{factorization} (ii), we get an open covering $B=\cup B_h$ of $B$ and morphisms $t_h\col B_h\ra\Ph^{ss}$ such that $t_h(b)$ is the configuration of lines associated to a principal part of $Y_b$ for $b\in B_h$. Consider the quotient morphism $q: p(\Gamma_\psi)\cap\Ph^{ss}\ra\J$. Now, $t_h(B_h\cap B_k)$ and $t_k(B_h\cap B_k)$ are congruent (in an algebraic way) modulo $SL(3),$ i.e. we can glue the morphisms $q\circ t_h$ to a morphism $B\ra\J$. By (\ref{orbit}), this morphism does not depend on the principal part of $Y_b$, hence we get the required functor transformation. Now, $P_g$ coarsely represents $\mathcal P_g$, thus we get a  morphism $P_g\ra\J$ which is $S_{m_g}$-invariant by construction, hence we get also a morphism $\beta\col P_g/S_{m_g}\ra\J$.
\end{proof}

Let us recall how the morphism $F_g\col\N\ra\B$ is defined in \cite{AL}. Let $(Y,p_1,\dots,p_{m_g})\in\N$.  The \emph{weighted dual tree} of $Y$ is the dual graph $\Gamma_Y$ of $Y$ and, for each vertex, the number of marked points contained in the corresponding component. If $\Gamma$ is a subset of $\Gamma_Y,$ let $wt(\Gamma)$ be the sum of the weights of the vertices contained in $\Gamma.$ We say that $Y$ has a \emph{central vertex} $v\in\Gamma_Y$ if $wt(\Gamma)<m_g/2$ for every connected subsets $\Gamma$ of $\Gamma_Y\backslash v$. The following is \cite[Lemma 3.2]{AL}.

\begin{Lem}\label{central-vertex}
Let $(Y,p_1,\dots,p_{m_g})\in\N$. Then $Y$ has a central vertex if and only if there are no edges $e$ of $\Gamma$ such that $wt(\Gamma_1)=wt(\Gamma_2)=m_g/2$, where $\Gamma_1,\Gamma_2$ are the connected subgraphs such that $\Gamma_1\cup\Gamma_2=\overline{\Gamma-e}$. There exists at most one central vertex.
\end{Lem}

A marked stable curve has a central vertex if it is not contained in the divisor 
$\Delta$ of $\N$ whose general point has two components containing $g+1$ marked points. If the central vertex exists, then $F_g(Y,p_1,\dots p_{2g+2})$ is obtained by contracting all the components of $Y$, which do not correspond to the central vertex. If $Y$ has no central vertex, it is easy to see that the edge  disconnecting $\Gamma_Y$ in two subgraphs with weights $g+1$ is unique.  We call it \emph{the central edge} of $\Gamma.$

\begin{Thm}\label{fact}
The variety $\J$ is a compactification of $H_g$ and the chain of maps 
$\N\stackrel{\beta_g}{\dashrightarrow}\J\stackrel{\alpha_g}{\dashrightarrow}\B$ gives a rational factorization of 
$F_g:\N\ra\B$.
 \end{Thm}

\begin{proof}
The morphism $\beta_g$ restricts to an injection $H_g\ra\J$ whose image is the subset $J_g$ of Theorem \ref{first}. The inverse is the morphism of Theorem \ref{first}.   

To prove the factorization, it is enough to show that, if $P$ is an irreducible principal part of a pointed stable curve $(Y,p_1,\dots,p_{m_g})$, then $P$ is the component corresponding to a central vertex of $\Gamma_Y$. First of all, assume that $Y$ has a (unique) central vertex $v$ and let $C\subset Y$ be the corresponding component. Assume that $P\ne C$ and let $v_P$ be the vertex of $P$. There is an edge $e$ of $\Gamma_Y$ such that, if $\Gamma_1,\Gamma_2$ are the connected subgraphs with $\Gamma_1\cup\Gamma_2=\overline{\Gamma_Y-e},$ then $v_P\in\Gamma_1$ and $v\in\Gamma_2$. Thus $wt(\Gamma_1)<m_g/2$ by definition of central vertex. Let $T$ be the conic twister such that $\omega_Y^\vee\otimes T$ has degree 2 on $P$. Now, $\psi(\phi_T(p_1),\dots,\phi_T(p_{m_g}))\in\Ph^{ss}$, by definition of principal part. Since $\phi_T$ contracts the connected subcurve of $Y$ corresponding to $\Gamma_2$ to a unique marking, we have $wt(\Gamma_2)<m_g/2$ by Lemma \ref{unstable} (i). Then $m_g=wt(\Gamma)=wt(\Gamma_1)+wt(\Gamma_2)<m_g,$ a contradiction. If $Y$ has no central edge, let $e$ be the central edge of $Y$. Set $\Gamma_1\cup\Gamma_2=\overline{\Gamma_Y-e}$ for connected graphs $\Gamma_1,\Gamma_2$. If $T$ is the twister of $P,$ then $\phi_T$ contracts either the component of $\Gamma_1$ or of $\Gamma_2$. Now, $\psi(\phi_T(p_1),\dots,\phi_T(p_{m_g}))\in\Ph^{ss}$ and $(\phi_T(p_1),\dots,\phi_T(p_{m_g}))$ has a marking of weight at least $g+1$, because $wt(\Gamma_1)=wt(\Gamma_2)=g+1$, contradicting Lemma \ref{unstable} (i).
\end{proof}

\subsection*{Acknowledgments} I wish to thank Lucia Caporaso and Edoardo Sernesi for useful conversations and helpful comments.

\end{document}